\documentclass[11pt]{amsart}
\usepackage{float}
\usepackage{amsfonts, amstext, amsmath, amsthm, amscd, amssymb, upgreek}
\usepackage{graphicx, color,  subfigure, wrapfig, overpic}
\usepackage[all,cmtip]{xy} 
\usepackage{tikz}
\usepackage[toc,page]{appendix}
\usepackage{rotating}
\let\oldmarginpar\marginpar
\renewcommand\marginpar[1]{\oldmarginpar[\raggedleft\footnotesize #1]%
{\raggedright\footnotesize #1}}

\AtBeginDocument{
   \def\MR#1{}
}
\newcommand{\Sp}{\rm{S}} 
\newcommand{\C}{\mathbb{C}}
\newcommand{\R}{\mathbb{R}}

\newcommand{\Z}{\mathbb{Z}}

\newcommand{\HH}{\mathbb{H}}

\theoremstyle{plain}
\newtheorem{theorem}{Theorem}[section]

\newtheorem*{namedtheorem}{\theoremname}
\newcommand{\theoremname}{testing}

\theoremstyle{definition}
\newtheorem{define}[theorem]{Definition}

\newtheorem{remark}[theorem]{Remark}

\title[ Compl\'{e}ment  to the Thurston 3D-Geometrization ] { Compl\'{e}ment to the Thurston 3D-Geometrization Picture}
\author[Alice Kwon and Dennis Sullivan]{Alice Kwon and Dennis Sullivan}

\begin{document}
\maketitle
\begin{abstract}
  Geometrization says `` any closed oriented three-manifold which  is prime (not a connected sum) carries  one of the eight Thurston geometries OR  it has  incompressible torus walls whose complementary components each carry one of four particular Thurston geometries" (see Introduction and Figure 1). \\
  \indent These geometric components have finite volume for the hyperbolic geometries (the H labeled vertices). They also have finite volume for each of the two  geometries appearing as Seifert fibrations (the S labeled vertices). The remaining pieces (the I labeled vertices) have Euclidean geometries of linear volume growth. \\
  \indent Then these vertex geometries are combined topologically to recover the original manifold. This, by cutting off the toroidal ends and then gluing the torus boundaries by affine mappings (indicated by the labeled edges in Figure 1). \\
  \indent The point of this work is to make the affine gluing respect an interpretation of the metric geometry in terms of a new notion of `` regional Lie  generated geometry". The vertex regions use four geometries in Lie form  combined in the overlap edge regions  via  affine geometry.\\
  \indent The Theorem solves, using Geometrization, a 45 year old question/approach to the Poincar\'{e} Conjecture. This was described in a '76  Princeton Math dept. preprint and finally documented in the 1983 reference by Thurston and the second author.
   
 \end{abstract}  

\begin{figure}[H]
    \includegraphics [width=5cm]{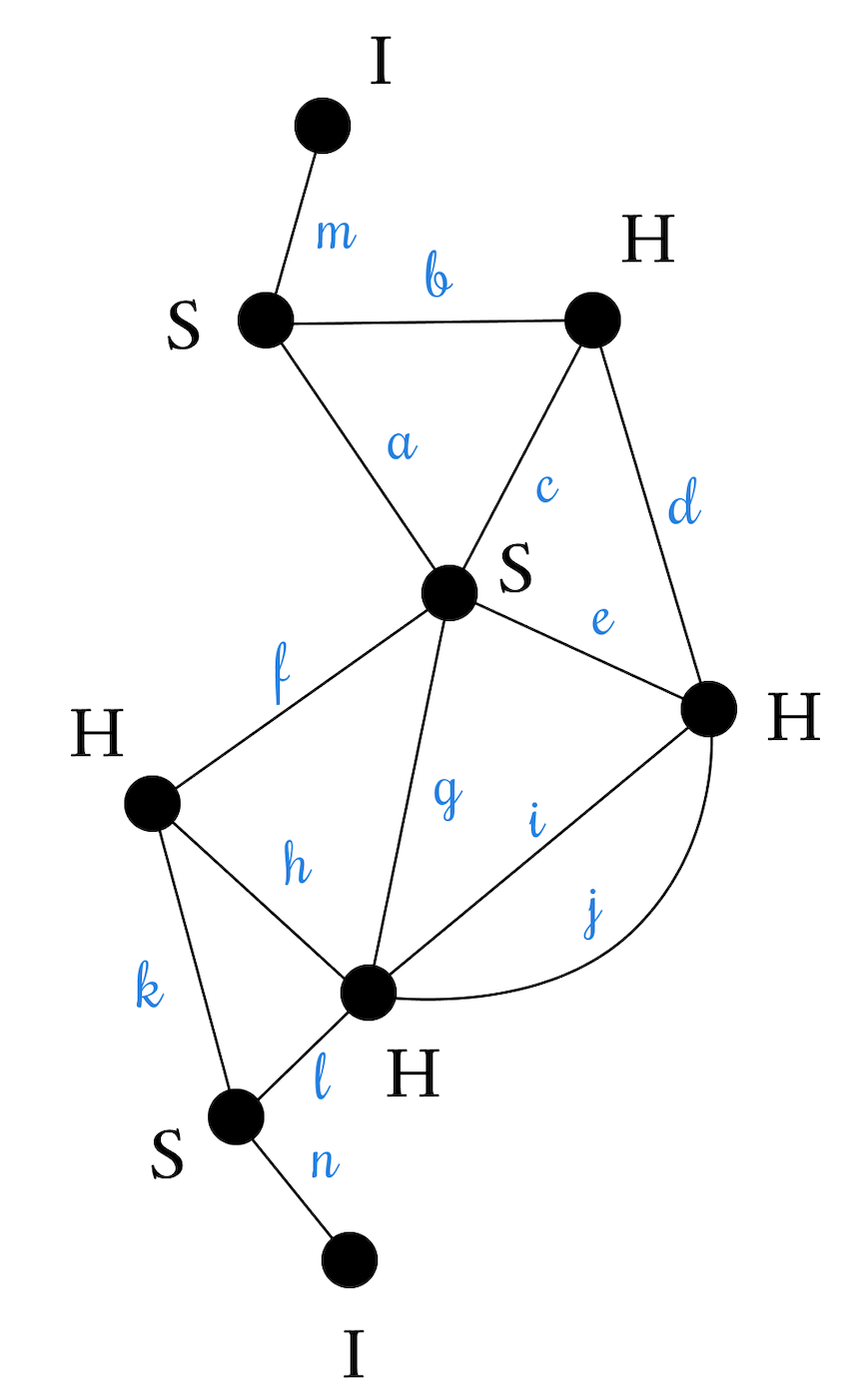}
    \caption{A closed oriented prime three-manifold made of Seifert fibrations (labeled S) and hyperbolic manifolds (labeled H)   and I-bundles over Klein bottles (labeled I) glued up by affine torus mappings (labeled $a,b,c,...$).}
    \label{fig:figureDennis3-1}
\end{figure} 
 \newpage
\section{Introduction }

 Certain {\bf prime} closed three-manifolds, which we call {\bf ``elemental" primes}, are by definition quotients by discrete  subgroups of isometries of one of the eight Thurston homogeneous Riemannian geometries. Recall that prime in 3D means {\bf not} a non-trivial connected sum. These eight 3D geometries naturally divide into three types. There are three of ``Hyperbolic type": hyperbolic three space, hyperbolic two space cross real line and the universal cover of  the  unit tangent bundle of  hyperbolic two space. There are three geometries of ``Affinely flat type": Euclidean three space, Heisenberg geometry, Solv geometry. Lastly, there are two of ``Elliptic type": the three-Sphere and the two-Sphere cross real line \cite{Scott}.\\
 \indent The elemental primes of the three ``hyperbolic types"  and of finite volume can be combined {\bf smoothly} according to any graphical pattern  with arbitrary torus gluing data. To an external half edge of the graph one can also attach a twisted I-bundle over a Klein bottle which will have a torus boundary. By these constructions all primes which are not elemental primes are obtained topologically \cite{BF}, \cite{P2}, \cite{Scott}. We call these manifolds \textbf{compound primes.}  The point of the paper is to construct these compound primes {\bf geometrically} using the six Lie groups corresponding to the six geometries of hyperbolic or affine type. Thus we build specific {\bf \textbf{Lie generated geometries}} on all of the compound prime closed oriented three-manifolds.  

\ \\
 \noindent {\bf Main Theorem.} Excepting those covered by $\Sp^3$ or $\Sp^2 \times \R$,
every closed oriented prime three-manifold $M$ has a canonical nonempty finite dimensional moduli space of Lie generated geometries defined by four Lie groups of dimensions (6,6,4,4) acting on the upper half space. 

\ \\
 \textbf{To be noted}: This paper takes advantage of a century of extremely great work using three additional points. The great work comprises the achievements of topologists, geometers and analysts: Poincar\'e's groupoid from classical multivalued functions defined by differential equations and his Conjecture... Seifert fibrations, Van Kampen-Milnor unique  connected sum decomposition of three-manifolds into primes, Stallings' fibration over the circle theorem, the Haken-Waldhausen incompressible submanifold techniques, the canonical JSJ decompositions, the Scott core... the Geometrization conjecture and proof by Thurston and Hamilton-Perelman. The additional points are: firstly, the very simple geometry of the very complicated irregular covers of the three-manifold corresponding to  the incompressible torus walls, secondly, the possibility to geometrically fit together metric Thurston pieces as Lie objects  using affine interpolation between the non compact ends of quotients of the four Lie groups when the four are engineered into one model space with the ends of linear form, and  thirdly, the notion of a Lie generated geometry which conceptually organizes the first two points. The Theorem solves a 45 year old question documented in \cite{DT}.
 
\section{Background and Definition of Lie generated geometry.}
 
 \textbf{Background}:\\
i) Recall \cite{DT}, if $G$ is a (finite dimensional) Lie group acting on a model  manifold $U$ real analytically, a $(G,U)$ structure on a manifold $M$ of the same dimension as $U$ is \textbf{determined} by an open cover of $M$ by coordinate charts in the model space $U$ whose partially defined transition mappings in $U$ extend globally and uniquely to $U$ being restrictions of elements of $G$ acting on $U$. The discussion that follows can also be applied to any group acting real analytically on the model space $U$. \\
ii) We will say that two $(G,U)$ structures agree if and only if they are equal, after forming the atlas of all charts of $M$ into $U$ which are $G$ compatible as mentioned above.\\
iii) Ehresmann (1950), cf.\cite{DT, ehressmann}, went one step further by passing to the set $\mathcal{G}$ of all germs of charts in $U$ at  all of the various points of $M$ of this  maximal atlas of charts.\\
iv) Each chart in the atlas defines a basic open set for a topology on  $\mathcal{G}$.  This topology defines   $\mathcal{G}$ as {\bf a sheaf of  germs of charts} over $M$  whose vertical fibers are  acted on freely and transitively by $G$.  The topology in the $G$ direction is discrete.\\
v) Each path in $M$ defines a $G$ equivariant map between the fibres over its endpoints. Namely each path lifts canonically for each germ over the starting point by analytic continuation of the germ of chart along the path. The germ of chart at the endpoint of the path only depends on the homotopy class of paths up to continuous deformation  fixing the endpoints.    \\
vi) The path components of $\mathcal{G}$ are regular or Galois  covering spaces of $M$ on which the multivalued analytic continuation maps into $U$ are single-valued. These maps on the  path components of $\mathcal{G}$ are referred to as the {\bf developing mappings} of the $(G,U)$ structure on $M$. These mappings are equivariant with respect to the holonomy of the Poincar\'e path groupoid. The latter is seen by noting the effect of continuation along paths.

\begin{define}{\bf Lie generated geometry}
\\
We consider several Lie groups $\{G\}$ acting real analytically on one model manifold $U$ and define, based on the background discussion, the notion of a Lie generated geometry $(\{G\},U)$ on $M$. Here is the definition.\\
 i) Suppose $M$ is covered by open connected regions $\{R\}$ and each is provided with a $(G,U)$ structure (as in the Background) from the collection of actions $(\{G\},U)$.\\
 ii) We can afford for the work here to restrict to cases where only pairs of regions can intersect. Thus, there are no triple intersections. \\
 iii) On each binary intersection $V$ one has by restriction, from the two regions which intersect in $V$, two (sheaves of germs of charts in $U$) structures, say $\mathcal{G}$ and $\mathcal{G}'$.  Each one consists of path components of germs which are regular covering spaces of $V$. \\
 iv) (key axiom) Assume there is a third ($H$,$U$) structure on $V$ whose germs of charts are contained in the intersection of  $\mathcal{G}$ and $\mathcal{G}'$, as subsets of all germs of charts into $U$. This means the $H$ action on $U$ lies in the intersection of the $G$ and the $G'$ actions on $U$. The ($H$,$U$) structure is termed a \textbf{reduction} of each of the two structures $\mathcal{G}$ and $\mathcal{G}'$.\\
 v) Once the regions are precisely specified, we say that two Lie generated geometries agree when the following construction produces equal results. Firstly, for each region form the set of germs of charts for that region. Secondly, consider for each binary intersection of regions the germs of charts in $U$ for the {\bf minimal} ($H$,$U$) structure  that  satisfies the essential axiom iv). This completes the definition of a Lie generated geometry.
 
\ \\
\textbf{ Derived Structure}: One may also  consider the derived or  associated (single group, $U$) structure on $M$ obtained by forming the group action generated by the set of Lie group actions on $U$. Making this structure canonical and concrete is done by saturating the germs at play in the original definition  by the action of the one big group action  which is generated.  Denote this action ($<\negmedspace \{G\} \negmedspace>$, $U$). This derived  action and structure is useful to check properties of the developing mappings and defining holonomy. This structure has less information, e.g. the finite dimensional moduli space of Lie generated structures   mentioned in the Theorem. The moduli space for this derived structure would tend to be infinite dimensional if one considered only $<\negmedspace \{G\} \negmedspace>$ directly. 
\end{define}


\ \\

\begin{remark} \label{rmk:two}
We will use, back and forth, the canonical real analytic identification between euclidean three space $\R^3 = \{x,y,z\}$ and upper half space $U = \{x,y,z|z>0\}$ given by the $\log$ and $\exp$ on the ($z$ coordinate)  $\times$ (the identity of the $x$ and $y$ coordinates). This transfers the Lie group denoted $[A]$ of dimension 6, acting affinely on $\R^3$ [which respects the $\R^2 \times \R^1$ product structure and is measure preserving on each factor] to a subgroup of {\bf affine  transformations} acting on $U$. This is one of the four Lie groups acting on $U$ in the Theorem. Two other Lie groups of the Theorem are ${\rm PSL}(2,\C)$ (of dimension 6) and ${\rm PSL}(2,\R) \times \R$ (of dimension 4). Finally, the fourth group of the Main Theorem is a $\Z$ cover of ${\rm PSL}(2,R)\times {\rm SO}(2)$ (of dimension 4).
\end{remark}

\indent {\it Proof of Main Theorem.}
 We first study the four Lie group actions on $U$ and the toroidal ends of various discrete subgroup quotients. The conclusion is that for these actions the $\Z + \Z$ subgroup corresponding to the possible toroidal ends are conjugate in each of the four cases to lattices in the group of horizontal translations of upper half space fixing infinity. We assume this for now. The proof, which is easy to find, except in one case, is contained in Discussion A below.
\ \\
   Step one: We consider, assuming this standard form, an elemental prime of finite volume or of linear growth with its $(G,U)$ structure, a toroidal end and its associated flat torus $T$  which is chosen to have normalized area, possible by Discussion A.
\begin{figure}[H]
       \includegraphics [width=12cm]{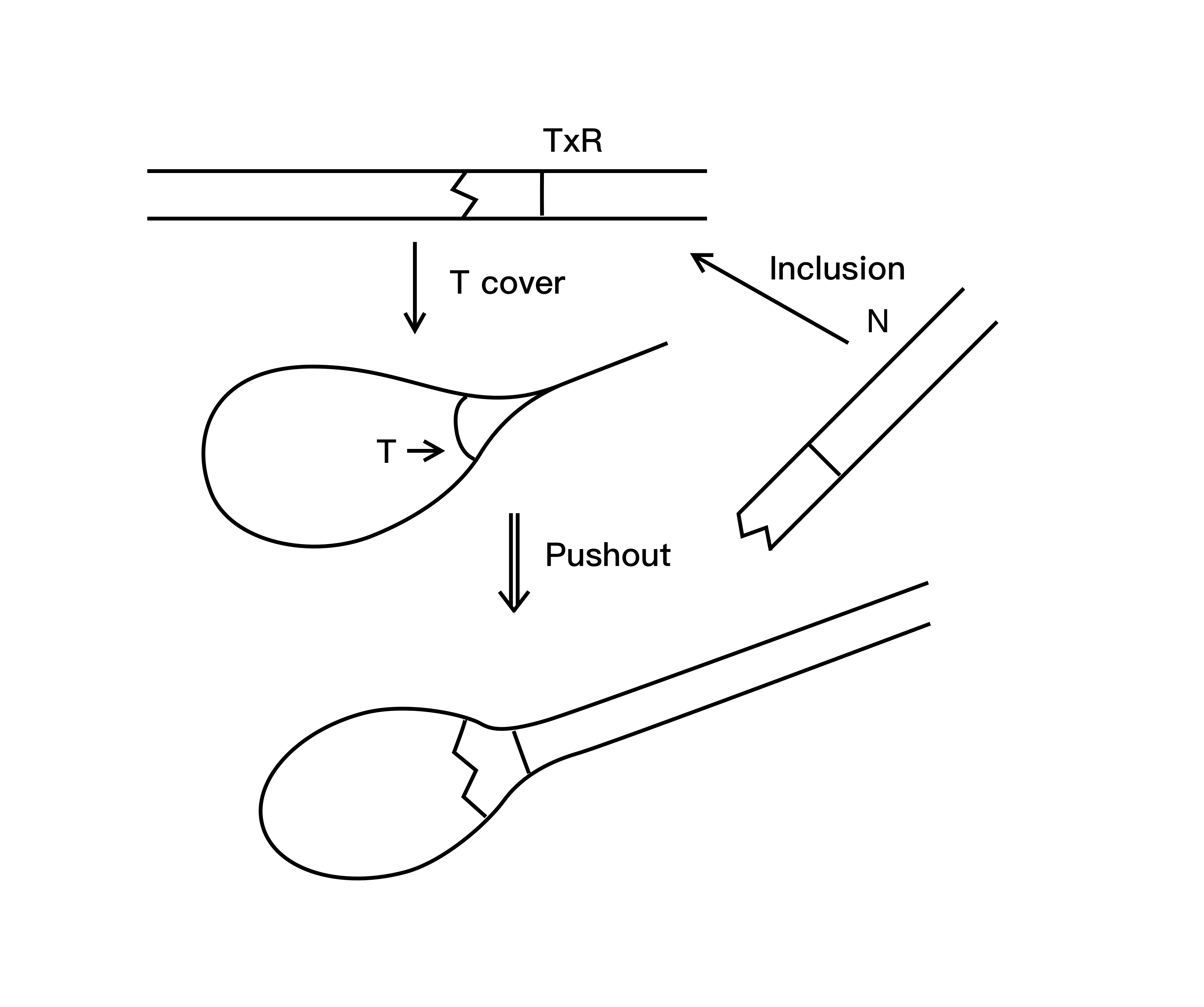}
   \caption{Attaching a flat cylinder to a finite volume Thurston geometry as a Lie generated geometry}   
   \label{fig:figureDennis1-1}
 \end{figure}     
   
   We form a Lie generated geometry using the idea of Figure \ref{fig:figureDennis1-1}. The $T$ cover arrow in Figure \ref{fig:figureDennis1-1} is the covering space associated to the $\Z + \Z$ subgroup of the manifold coming from the toroidal end with basepoint on $T$.  Remarkably, the total space of this very complicated irregular cover is a simple flat manifold, being $U$ mod $\Z + \Z$, to be sure after the third coordinate is stretched out by the $\log$ as discussed in Remark \ref{rmk:two}. This flat manifold has the very simple form $T\times \R$. We consider the neighborhood $N$ of the right hand end with boundary as indicated slightly to the left of $T \times$ a point. This $N$ is the component of the pre-image of the end in the manifold beyond $T$ containing the canonical base point above the base point.\\
  \indent The slanting arrow is the inclusion of a copy of $N$ into $T \times \R$. The pushout of this diagram  of two arrows (the $T$ cover arrow and the inclusion arrow) is diffeomorphic to the original manifold. The push out manifold  inherits  the $(G,U)$ structure of the original manifold on the region $\R$ which  consists of points to the left of $T$ plus the reflection of the indicated collar in the upper $T \times \R$ to the right of $T$. The region $N$ inherits by inclusion the flat affine structure .\\
   \indent Denote by $V$ the double collar  of $T$. Then $V$ gets two structures by restriction of those from $R$ and from $N$. Each of these contains the $(\Z+ \Z,U)$ structure  on $V$ defined by the lattice determining the flat torus $T$. This creates a Lie generated structure on $M$ which attaches a half infinite flat cylindrical torus end (of linear growth) in place of the original finite volume end. This was our goal in this first step.
   
   \begin{figure}[H]
       \includegraphics [width=13cm]{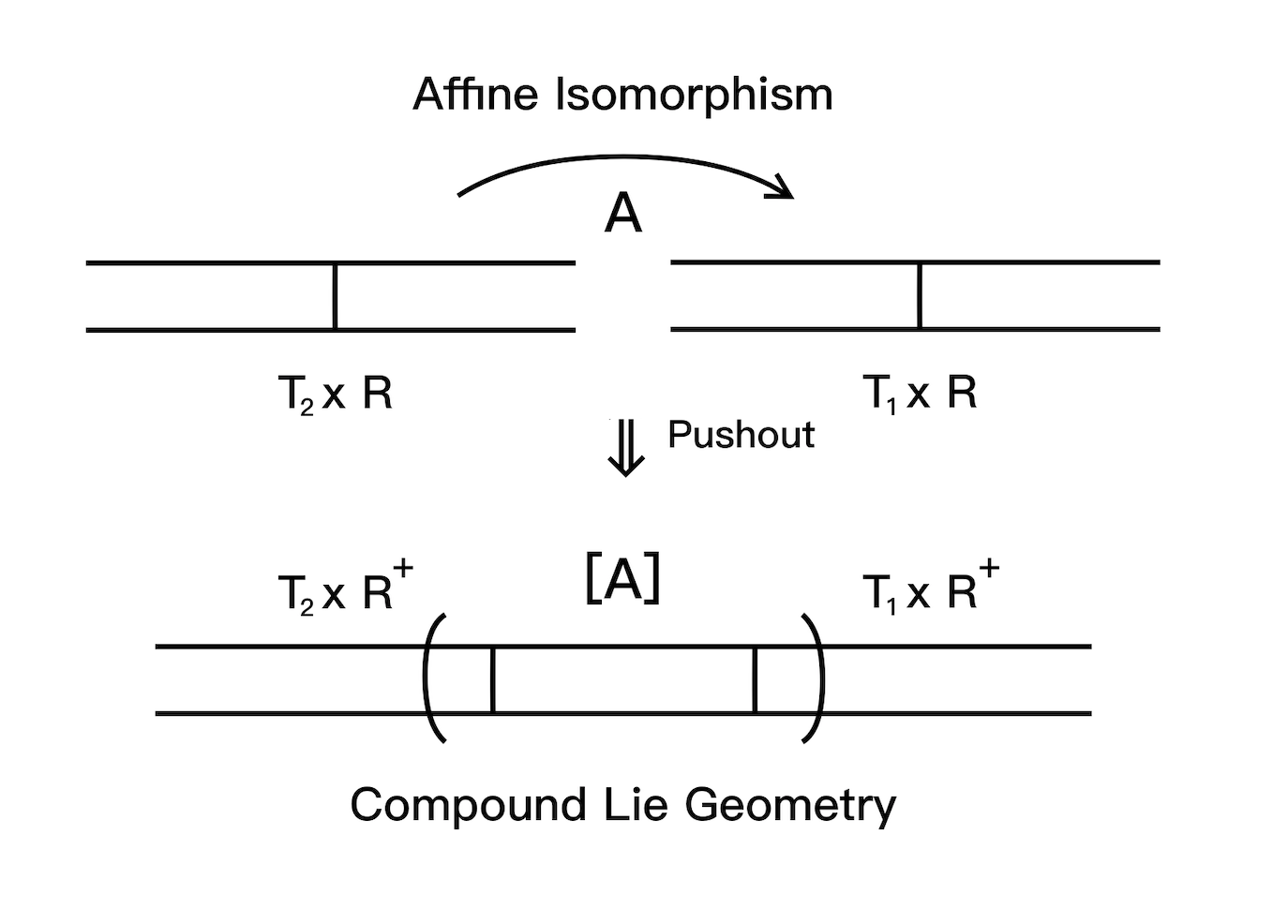}
     \caption{The affine cylinder incorporating the torus gluing data: the two flat torii and a class of diffeomorphisms between them.}
     \label{fig:figureDennis2-1}
 \end{figure} 
 
 \ \\
   \indent Step two: Consider two flat tori of normalized area, $T_1$ and $T_2$, and a homotopy class of  diffeomorphisms between them. This class is represented  by an area preserving affine mapping denoted $A$.\\
\indent We consider, looking at Figure \ref{fig:figureDennis2-1}, gluing by $A$ the boundary of the left half of $T_2 \times \R$ to the right half of $T_1 \times \R$.  Actually we glue a symmetric neighborhood of $T_2$ to a symmetric neighborhood of $T_1$ by the restriction of $A \times$ reflection in the third coordinate. The result is pictured in the lower component of Figure \ref{fig:figureDennis2-1}. It is a Lie generated geometry with three regions, left, right and middle, denoted $L$, $R$ and $M$, see next paragraph.\\
   \indent The structure on this manifold is generated by the two lattices and the affine conjugacy flip between them, to be explained now. \\
   \indent The entire group $H$ which these lattices and the flip $\times A$ generate is contained in the affine group $[A]$ of dimension 6 referred to in Remark \ref{rmk:two} and Discussion A). From the point of view of $[A]$ the Lie generated geometry constructed can be seen as follows: firstly provide a standard Torus $\times \R$ with the affine structure corresponding to the ($[A]$, $U$) action, secondly divide the Torus $\times \R$ into three regions $L$,$R$ and $M$  as indicated above and thirdly consider reductions of the $[A]$ structure to the two different lattices respectively at each end $L$ and $R$. In the middle $M$ choose the reduction to the structure corresponding to the subgroup $H$ of $[A]$ generated by these two lattices and the flip conjugacy $A \times$ reflection.\\
   \indent One obtains by these reductions a Lie generated geometry on the Torus $\times \R$ with given toroidal cylinders on the two ends glued by any prescribed isotopy class of diffeomorphisms between the tori. This was our goal in step 2).
   
   \ \\
   \indent Step three: To treat interval bundles over a Klein bottle with a  torus  boundary, we proceed as follows: Firstly form a standard $T \times \R$.  Secondly, construct an orientation preserving fixed point free affine involution of $T \times \R$, that flips the ends and that keeps one $T$ invariant, where it reverses the orientation.  The quotient by this fixed point free involution will be an orientable three-manifold which is a twisted $\R$-bundle over a Klein bottle with one end and that end is geometrically the same as one end of  $T \times \R$. This is the goal of step 3.
   
   \ \\
   \indent Step four: Now we finish the proof of the Main Theorem assuming the standard form demonstrated in Discussion A.  Take any connected finite graph. Decorate some vertices with finite volume Thurston elemental primes. The valence of each such vertex corresponds to the number of toroidal ends of the finite volume manifold which decorate the half edges at the vertex. The graph may also have some univalent vertices where  twisted $\R$ bundles over Klein bottles are assigned. Finally, each edge of the graph is assigned an isotopy class of diffeomorphisms between the torus boundary components corresponding to the half edges.\\ 
   \indent From the geometrization Theorem \cite{BF} one knows an oriented prime either carries a Thurston geometry or it admits a graphical toroidal decomposition with the properties we have just described. In each case depending on finitely many parameters.  
   
   To each full edge we take note of the two flat tori and the  diffeomorphism class of gluing data and form the Affine object of step two and Figure \ref{fig:figureDennis2-1}, with which we decorate the edge.

  Now there is for each half edge of the entire graph an identification to perform.
  But now the flat tori match exactly and the identification is by the identity on these two copies of the same torus. We form these identifications by sliding foreshortened half cylinders together rigidly preserving the geometry. How much sliding one does adds to the list of finite parameters.\\
  \indent This completes the proof of the Main Theorem after checking Discussion A.

  \section{Discussion A}
   \indent  i) For the affine group of the four Lie groups of the Main Theorem \ mentioned in Remark \ref{rmk:two}, after identifying $\R^3$ with upper half space $U$ by the $\exp$ of the third coordinate, one sees immediately the $\Z + \Z$ subgroups are lattices of horizontal translations fixing infinity. Call these of standard form, recalling they depend on a parameter, the classical $j$-invariant.\\ 
   
   \indent ii) For a cuspidal end of an $\HH^3$ finite volume manifold which when represented as a $U$ quotient with one lift of the cusp at infinity in $U$ one sees immediately the $\Z + \Z$ subgroup of the cusp is of standard form. Now the lattice is  somewhat mysterious and   general although arithmetic in an  appropriate sense. Note that the induced area of the torus takes on at least an interval of values $(0,c]$ for $c$  some universal constant. This is convenient and  means that we can choose the height of the lattice  in $U$ so that its area is normalized.\\
   
 \indent  iii) For the $\HH^2 \times \R$ geometry, which acts on $U$ = upper half plane $\times \R$ this check is again immediate that the $\Z + \Z$ subgroup is of standard form. The shape of the lattice is now rectangular, with one side length taking on  at least an interval of values $(0,c']$ and the other of any chosen normalized length. Thus again the area of the torus takes on a universal interval of values abutting zero which allows  us to normalize the area.\\
 
 \indent This  finishes the proof in three of the Lie group cases that the action on $U$ of any toroidal end group $\Z + \Z$ is of standard form, and ends the first part of Discussion A. For the  fourth and last case  some  additional work is needed.
   
 \ \\
  \indent iv) The orientation preserving (both fibre and base) isometry group $G'$ of the universal cover $T'$ of the unit tangent bundle $T$ of the hyperbolic plane is a particular $\Z$ cover of the group $G$  of orientation preserving isometries of the upper half plane acting on the unit tangent bundle $T$ cartesian product with the group ${\rm SO}(2)$ acting by  the frame bundle action on $T$, namely rotating the fibres  by constant angle \cite{BF}.  The $\Z$ cover $G'$ is defined by the kernel of the homomorphism of addition $\Z + \Z \to \Z$  where $\Z + \Z$ is identified 
  to the fundamental group of ${\rm PSL}(2,\R) \times {\rm SO}(2)$  using standard generators.

    Here is the  analysis of the action of the $\Z + \Z$ associated to a cuspidal end.
    Consider the cross-section  of the unit tangent circle bundle of $\HH^2$ defined by the unit vertical vector field in the upper half plane model. This  cross section  is preserved by the subgroup of isometries which fix infinity. This subgroup is generated by horizontal translations of $\HH^2$ and by dilations from any point on the real line.
    
    The $\Z$ corresponding to the fibre circle of the torus is translation in the $\R$ factor of $\HH^2 \times \R$, which is $U$ the upper half space. The horocycle circle gives the other $\Z$, which acts by the horizontal translation group in $\HH^2$. Fortunately this  translation amount is constant as we move out the cuspidal end and in this coordinate parametrization, the cuspidal end looks like a regular product toroidal cylinder.
    
     End of Discussion A and the proof of the Main Theorem. \qed

.
\bibliographystyle{plain}
\bibliography{references-ak}

\end{document}